\documentclass[a4paper]{article}
\usepackage[english]{babel}
\usepackage{amsmath,amssymb,amsfonts,amsthm}
\usepackage{mathrsfs}
\usepackage{graphicx}
\usepackage[numbers,sort&compress]{natbib}

\usepackage[bookmarks=false,%
pdfstartview=FitH,linkbordercolor={0.5 1 1},%
citebordercolor={0.5 1 0.5},unicode,%
pagebackref,hyperindex]{hyperref}%

\newtheorem{theorem}{Theorem}
\newtheorem{lemma}{Lemma}

\newtheorem{example}{Example}

\newcommand{\conv}{\mathop{\mathrm{conv}}}
\newcommand{\tr}{\mathop{\mathrm{tr}}}
\newcommand{\setA}{\mathscr{A}}
\newcommand{\setP}{\mathscr{P}}
\newcommand{\setV}{\mathscr{V}}
\newcommand{\setAinf}{\mathscr{A}_{\infty}}
\newcommand{\len}{\mathop{\mathrm{len}}}
\newcommand{\diag}{\mathop{\mathrm{diag}}}
\newcommand{\dist}{\mathop{\mathrm{dist}}}
\DeclareMathOperator*{\infx}{\vphantom{sup}inf}

\begin{document}
\date{}

\title{On explicit a priori estimates\\
of the joint spectral radius\\ by the generalized Gelfand
formula\thanks{This work was supported by the Russian
Foundation for Basic Research, project no.
06-01-00256.}%
~\thanks{Dedicated to memory of Bernd Aulbach}}

\author{Victor Kozyakin\\[5mm]
Institute for Information Transmission
Problems\\ Russian Academy of Sciences\\ Bolshoj Karetny lane 19,
Moscow 127994 GSP-4, Russia}

\maketitle

\begin{abstract}
In various problems of control theory, non-autonomous and
multivalued dynamical systems, wavelet theory and other fields
of mathematics information about the rate of growth of matrix
products with factors taken from some matrix set plays a key
role. One of the most prominent quantities characterizing the
exponential rate of growth of matrix products is the so-called
joint or generalized spectral radius. In the work some explicit
a priori estimates for the joint spectral radius with the help
of the generalized Gelfand formula are obtained. These
estimates are based on the notion of the measure of
irreducibility (quasi-controllability) of matrix sets proposed
previously by A. Pokrovskii and the author.

\medskip\noindent PACS number 02.10.Ud; 02.10.Yn


\medskip

\end{abstract}

\section{Introduction}\label{S-intro}
In various problems of control theory \cite{BrayTong:TCS80},
non-autonomous and multivalued dynamical systems
\cite{Bar:AIT88-2:e,AulSieg:01,AulSieg:02}, wavelet theory
\cite{DaubLag:SIAMMAN92,DaubLag:LAA92,ColHeil:IEEETIT92} and
other fields of mathematics information about the rate of
growth of matrix products with factors taken from some matrix
set plays a key role. One of the most prominent values
characterizing the exponential rate of growth of matrix
products is the so-called joint or generalized spectral radius.
The aim of the paper is to obtain efficient explicit estimates
of the joint spectral radius.

Let $\setA=\{A_{1},\ldots,A_{r}\}$ be a set of real $m\times m$
matrices. As usual, for $n\ge1$, denote by $\setA^{n}$ the set
of all $n$-products of matrices from $\setA$; $\setA^{0}=I$.

Let $\|\cdot\|$ be a fixed but otherwise arbitrary norm in
${\mathbb{R}}^{d}$. Then the limit
\begin{equation}\label{E-JSRad}
\hat{\rho}({\setA})=
\limsup_{n\to\infty}\|\setA^{n}\|^{1/n},
\end{equation}
where
\[
\|\setA^{n}\|=\max_{A\in\setA^{n}}\|A\|=\max_{A_{i_{j}}\in\setA}
\|A_{i_{n}}\cdots
A_{i_{2}}A_{i_{1}}\|,
\]
is called \emph{the joint spectral radius} of the matrix set
$\setA$ \cite{RotaStr:IM60}. In fact, the limit in
(\ref{E-JSRad}) does not depend on the norm $\|\cdot\|$.
Moreover, for any $n\ge 1$ the estimates
$\hat{\rho}({\setA})\le \|\setA^{n}\|^{1/n}$ hold
\cite{RotaStr:IM60}, and therefore the joint spectral radius
can be defined also by the following formula:
\begin{equation}\label{E-JSR}
\hat{\rho}({\setA})=
\inf_{n\ge 1}\|\setA^{n}\|^{1/n}.
\end{equation}

If the matrix set $\setA$ consists of a single matrix then
(\ref{E-JSRad}) turns into the known Gelfand formula for the
spectral radius of a linear operator. By this reason sometimes
(\ref{E-JSRad}) is called the generalized Gelfand formula
\cite{ShihWP:LAA97}.

For each $n\ge1$ it can be defined also the quantity
\[
\bar{\rho}_{n}(\setA)=\max_{A_{i_{j}}\in\setA}
\rho(A_{i_{n}}\cdots
A_{i_{2}}A_{i_{1}}),
\]
where maximum is taken over all possible $n$-products of the
matrices from the set $\setA$, and $\rho(\cdot)$ denotes the
spectral radius of the corresponding matrix, that is the
maximum of modules of its eigenvalues. In these notations the
limit
\begin{equation}\label{E-GSRad}
\bar{\rho}({\setA})=
\limsup_{n\to\infty}\left(\bar{\rho}_{n}(\setA)\right)^{1/n}
\end{equation}
is called \emph{the generalized spectral radius} of the matrix
set $\setA$ \cite{DaubLag:LAA92,DaubLag:LAA01}. Moreover, for
any $n\ge 1$ the estimates $\bar{\rho}({\setA})\ge
\left(\bar{\rho}_{n}(\setA))\right)^{1/n}$ hold, and hence the
generalized spectral radius can be defined also as follows:
\begin{equation}\label{E-GSR}
\bar{\rho}({\setA})=
\sup_{n\ge 1}\left(\bar{\rho}_{n}(\setA)\right)^{1/n}.
\end{equation}

As is shown in \cite[Thm. 2]{BerWang:LAA92}, see also
\cite{Els:LAA95,ShihWP:LAA97,Shih:LAA99}, the quantities
$\bar{\rho}({\setA})$ and $\hat{\rho}({\setA})$ coincide with
each other if the matrix set $\setA$ is bounded. This allows to
speak simply about \emph{spectral radius} of $\setA$, which
will be denoted in what follows as $\rho({\setA})$
($=\bar{\rho}({\setA})=\hat{\rho}({\setA})$).

In view of (\ref{E-JSR}), (\ref{E-GSR}) for any $n$ the
quantities $\left(\bar{\rho}_{n}(\setA)\right)^{1/n}$ and
$\|\setA^{n}\|^{1/n}$ form lower and upper bounds,
respectively, for the spectral radius of $\setA$:
\begin{equation}\label{Eq-sprad}
\left(\bar{\rho}_{n}(\setA)\right)^{1/n}\le
\rho({\setA})\le
\|\setA^{n}\|^{1/n}.
\end{equation}
This last formula may serve as the basis for a posteriori
estimating the accuracy of computation of $\rho({\setA})$. The
first algorithms of a kind in the context of control theory
problems have been suggested in \cite{BrayTong:TCS80}, for
linear inclusions in \cite{Bar:AIT88-2:e}, and for problems of
wavelet theory in
\cite{DaubLag:SIAMMAN92,DaubLag:LAA92,ColHeil:IEEETIT92}. Later
the computational efficiency of these algorithms was
essentially improved in \cite{Grip:LAA96,Maesumi:LAA96}.
Unfortunately, the common feature of all such algorithms is
that they do not specify the amount of computational steps
required to achieve desired accuracy of approximation of
$\rho({\setA})$.

In \cite{RotaStr:IM60,Els:LAA95} it was proved that the
spectral radius of the matrix set $\setA$ can be determined by
the equality
\begin{equation}\label{E-inf}
\rho({\setA})=\inf_{\|\cdot\|}\|\setA\|,
\end{equation}
where infimum is taken over all norms in ${\mathbb{R}}^{d}$.
For irreducible matrix sets\footnote{A matrix set $\setA$ is
called \emph{irreducible}, if the matrices from $\setA$ have no
common invariant subspaces except $\{0\}$ and
${\mathbb{R}}^{m}$.} $\setA$ infimum in (\ref{E-inf}) is
attained, and for such matrix sets there are norms $\|\cdot\|$
in ${\mathbb{R}}^{d}$, called \emph{extremal norms}, for which
\begin{equation}\label{E-extnorm}
    \|\setA\|\le\rho({\setA}).
\end{equation}

In various situations it is important to know the conditions
under which $\rho({\setA})>0$. To answer this question, in
general case one can use, for example, Theorem A from
\cite{Bochi:LAA03}, which claims the existence of a constant
$C_{d} > 1$ depending only on dimension of a space such that
for any bounded matrix set $\setA$ and any norm $\|\cdot\|$ in
$\mathbb{R}^{d}$ the following inequality holds:
\begin{equation}\label{E-Bineq}
\|\setA^{d}\|\le C_{d}\,\rho(\setA)\|\setA\|^{d-1}.
\end{equation}
From here it follows that the equality $\rho(\setA)=0$ implies the equality
$\|\setA^{d}\|=0$, and then also the equality $\setA^{d}=\{0\}$. In virtue of
(\ref{E-JSR}) the converse is also valid: $\setA^{d}=\{0\}$ implies
$\rho(\setA)=0$. So, theoretically verification of the condition $\rho(\setA)=0$
may be fulfilled in a finite number of steps: it suffices only to check that all
$d$-products of matrices from $\setA$ vanish. Of course this remark is hardly
suitable in practice since even for moderate values of $d = 3, 4$, $r= 5, 6$ the
computational burden of calculations becomes too high.

Remark that by (\ref{E-extnorm}) $\rho({\setA})>0$ for irreducible matrix sets
$\setA$.

Some works suggest another formulas to compute $\rho({\setA})$. So, in
\cite{ChenZhou:LAA00} it is shown that
\begin{equation}\label{E-tr}
\rho({\setA})=
\limsup_{n\to\infty}\max_{A_{i_{j}}\in\setA}
\left|\tr(A_{i_{n}}\cdots
A_{i_{2}}A_{i_{1}})\right|^{1/n},
\end{equation}
where as usual $\tr(\cdot)$ denotes the trace of a matrix. In
\cite{Maesumi:CDC05,Prot:FU98,Prot:CDC05-1}
$L^{p}$-ge\-ne\-ral\-iz\-ations of the formulas (\ref{E-JSR}),
(\ref{E-GSR}), (\ref{E-tr}) and (\ref{E-inf}) are proposed for
computation of the spectral radius of matrix sets. Algorithms
for computation of $\rho({\setA})$ based on the relation
(\ref{E-inf}) are considered, e.g., in
\cite{GugZen:CDC05,GugZen:LAA08,Maesumi:CDC05}. In
\cite{ParJdb:LAA08} it was noted that the norm in the
definition of the joint spectral radius can be replaced by  a
positive homogeneous polynomial of even degree. By developing
further this idea one can replace the norm in (\ref{E-JSRad})
by an arbitrary homogeneous  function strictly positive outside
of zero. Namely, let $\nu(x)$ be a strictly positive for
$x\neq0$ homogeneous function with degree of homogeneity
$\varkappa>0$, that is $\nu(tx)\equiv t^{\varkappa}\nu(x)$ for
any $t>0$. Then by introducing for an arbitrary matrix $A$ the
notation
\[
\nu(A)=\sup_{x\neq0}\frac{\nu(Ax)}{\nu(x)},
\]
one can easily get the following generalization of formula (\ref{E-JSRad}):
\begin{equation}\label{E-polyGSR}
    \rho({\setA})=
\limsup_{n\to\infty}\left(\nu_{n}({\setA})\right)^{1/\varkappa n},
\end{equation}
where
\[
\nu_{n}({\setA})=\max_{A_{i_{j}}\in\setA}
\nu(A_{i_{n}}\cdots
A_{i_{2}}A_{i_{1}}).
\]
As is shown in \cite{ParJdb:LAA08}, in a number of situations (\ref{E-polyGSR})
gives a better approximation to $\rho({\setA})$ because the collection of
positive homogeneous functions is richer than the collection of norms. In
particular, the Lebesgue sets\footnote{The Lebesgue set of a function $\nu(x)$
is the set
$\left\{x:~\nu(x)\le c\right\}$ for some $c$.}
of positive homogeneous functions may be non-convex in contrast to the Lebesgue
sets of norms.

In \cite{BN:SIAMJMA05} it is established that in the case when all matrices from
$\setA$ have non-negative entries the following inequalities hold
\begin{equation}\label{E-Bar}
    \frac{1}{r^{1/n}}\rho^{1/n}(A_{1}^{\otimes n}+\dots+A_{r}^{\otimes n})
    \le \rho({\setA})
    \le \rho^{1/n}(A_{1}^{\otimes n}+\dots+A_{r}^{\otimes n}),
\end{equation}
where $A^{\otimes n}$ denotes the $n$-fold Kronecker (tensor) product of the
matrix $A$ with itself. Here the fact that the right-hand and left-hand sides of
the inequalities (\ref{E-Bar}) do not contain mixed products of the matrices
from $\setA$  looks somewhat surprising. Theoretically, the inequalities
(\ref{E-Bar}) allow to compute $\rho({\setA})$ with any desirable accuracy.
However,
dimension of the matrix $A_{1}^{\otimes n}+\dots+A_{r}^{\otimes n}$ increases in
$n$ so rapidly that even for moderate values of $d = 3, 4$, $r = 5, 6$
computations become practically impossible. In the general case of arbitrary
matrix sets $\setA$ a bit more complicated analog of formula (\ref{E-Bar})
is also valid \cite{BN:SIAMJMA05}.

In investigation of properties of the spectral radius of matrix sets some
implicit definitions of the joint (generalized) spectral radius play an
important role. Let the matrix set $\setA$ be irreducible. Then
\cite{Bar:AIT88-2:e}
the value of $\rho$ equals to $\rho({\setA})$ if and only if if there is a norm
$\|\cdot\|$ in
${\mathbb{R}}^{m}$ such that
\begin{equation}\label{E-barnorm}
    \rho\|x\|\equiv\max_{A_{i}\in\setA}\|A_{i}x\|.
\end{equation}
The norm satisfying (\ref{E-barnorm}) is called
the \emph{Barabanov norm}. Similarly \cite[Thm. 
3.3]{Prot:FPM96:e}, \cite{Prot:FU98}, the quantity $\rho$ equals to
$\rho({\setA})$ if and only if for some central-symmetric convex
body\footnote{The set is called the \emph{body} if it contains at least one
interior point.} $S$ the following equality holds
\begin{equation}\label{E-protset}
    \rho S =\conv\left(\bigcup_{i=1}^{r}A_{i}S\right),
\end{equation}
where $\conv(\cdot)$ stands for the convex hull of a set. As is noted by V.
Protasov in \cite{Prot:FPM96:e}, the relation
(\ref{E-protset}) was proved by A. Dranishnikov and S. Konyagin, so it is
natural to call the central-symmetric set $S$ the
\emph{Dranishnikov-Konyagin-Protasov set}.
The set $S$ can be treated as the unit ball of some norm $\|\cdot\|$ in
${\mathbb{R}}^{d}$ (recently this norm is usually called the \emph{Protasov
norm}). As Barabanov norms as Protasov norms are the extremal norms, that is
they satisfy the inequality (\ref{E-extnorm}). In
\cite{PWB:CDC05,Wirth:CDC05,PW:LAA08} it is shown that  Barabanov and Protasov
norms are dual to each other.

Remark that formulas (\ref{E-inf}), (\ref{E-barnorm}) and (\ref{E-protset})
define the joint or generalized spectral radius for a matrix set in an
apparently computationally nonconstructive manner. In spite of that, namely such
formulas underlie quite a number of theoretical constructions (see, e.g.,
\cite{Koz:CDC05:e,Koz:INFOPROC06:e,Wirth:CDC05,Wirth:LAA02,ParJdb:LAA08,Bar:CDC05})
and algorithms \cite{Prot:CDC05-1} for computation of $\rho({\setA})$.

In \cite{ParJdb:LAA08} it was noted that the quantity $\rho$ is an upper bound
for $\rho({\setA})$ if for some strictly positive homogeneous polynomial $p(x)$
of degree $2d$ the following inequality holds:
\[
\max_{A_{i}\in\setA}p(A_{i}x)\le\rho^{2d}
p(x),\qquad\forall~x\in{\mathbb{R}}^{d}.
\]
By using this observation, in \cite{ParJdb:LAA08} algorithms for evaluating of
the spectral radius $\rho({\setA})$ were developed which are computationally no
less efficient and accurate than those proposed in \cite{BN:SIAMJMA05}.

As is shown in \cite[Lem. 2.3]{Wirth:IJRNC98}, \cite[Lem.
6.5]{Wirth:LAA02}, \cite[Sec. 5.2]{Prot:FPM96:e}, \cite[Thm.
3]{Prot:CDC05-1}, see also the survey in \cite{Theys:PhD05}, for irreducible
matrix sets $\setA$ the following inequalities are valid:
\begin{equation}\label{E-wirth-irr}
\gamma^{1/n}\|\setA^{n}\|^{1/n}\le \rho({\setA})\le\|\setA^{n}\|^{1/n}
\end{equation}
for a suitable constant $\gamma\in(0,1)$. Some computable estimates of the constant $\gamma$ are obtained in \cite[Sec. 8]{Prot:FPM96:e}.

When the matrix set $\setA$ is not irreducible the situation is more complicated. In this case, by the Bochi inequality (\ref{E-Bineq}), $\rho(\setA)$ may vanish, which happens if and only if $\setA^{d}=\{0\}$. But if $\rho(\setA)\neq0$ then \cite[Lem. 2.3]{Wirth:IJRNC98}
\begin{equation}\label{E-wirth-gen}
    \gamma_{*}^{(1+\ln
    n)/n}\|\setA^{n}\|^{1/n}\le\rho({\setA})\le\|\setA^{n}\|^{1/n}
\end{equation}
with some constant $\gamma_{*}\in(0,1)$, which is weaker than
(\ref{E-wirth-irr}). Unfortunately, \cite{Wirth:IJRNC98,Wirth:LAA02} contain
neither exact values nor at least effectively computable estimates for
$\gamma_{*}$.

The work is organized as follows. In Introduction we have presented a concise
survey of publications related to the problem of evaluation of the joint
(generalized) spectral radius. In Section~\ref{S-protf} the technique by V.
Protasov to get potentially computable estimates of the joint spectral radius is
described. In Theorem~\ref{T-main} from Section~\ref{S-irred} we obtain new a
priori estimates of the joint spectral radius based on the notion of the measure
of irreducibility (quasi-controllability) of matrix sets proposed previously by
A. Pokrovskii and the author. Section~\ref{PRT11} is devoted to the proof of
Theorem~\ref{T-main}. At last, in Section~\ref{S-ex} we cite examples from
\cite{KozPok:TRANS97}, which demonstrate how the value of the measure of
irreducibility for some matrix sets can be evaluated.

\section{Protasov Estimates}\label{S-protf}

Throughout the paper $\setA=\{A_{1},\ldots,A_{r}\}$ is an irreducible set of real $d\times d$ matrices. Provided that $\rho({\setA})=1$, one can derive from \cite{Prot:FPM96:e} the following estimate for the constant $\gamma$ in (\ref{E-wirth-irr}):
\begin{equation}\label{E-formprot}
\gamma\ge\frac{p_{1}({\setA})\cdots p_{d-1}({\setA})}{\left(1+\|\setA\|\right)^{d-1}},
\end{equation}
in which the quantities $p_{1}({\setA}),\dots, p_{d-1}({\setA})$
are determined by the equalities
\[
p_{k}({\setA})=\infx_{\substack{L\subset\mathbb{R}^{m}\\ \vphantom{\|}\dim L=k}}
\sup_{\substack{x\in L\\ \|x\|=1}}
\max_{A_{i}\in\setA}\dist(A_{i}x,L),\quad k=1,2,\dots,d-1,
\]
where the external infimum is taken over all subspaces
$L\subset\mathbb{R}^{d}$ of dimension $k$, $\|\cdot\|$ is the Euclidean norm\footnote{Formally, the norm $\|\cdot\|$ in \cite{Prot:FPM96:e} is not assumed to be Euclidean. However, to estimate the quantity $\widehat{H}=\gamma^{-1}$ in \cite[Thm. 8.2]{Prot:FPM96:e} the author uses the notion of perpendicularity of subspaces and evaluates the area of a triangle with sides measured by the norm $\|\cdot\|$, which implicitly mean that the norm $\|\cdot\|$ is Euclidean.}  in $\mathbb{R}^{d}$, and $\dist(A_{i}x,L)$
denotes the distance from the point $A_{i}x$ to the subspace  $L$, that is $\dist(A_{i}x,L)=\inf_{y\in L}\|A_{i}x-y\|$.

As is noted in \cite{Prot:FPM96:e}, all the quantities
$p_{1}({\setA}),\dots, p_{d-1}({\setA})$ are strictly positive since due to irreducibility of $\setA$ the matrices  $A_{1},\ldots,A_{r}$ have no non-trivial invariant spaces. Thus, the quantity $\gamma$ characterizes, in a sense, a `degree of irreducibility' of the matrix set $\setA$.

If $\rho({\setA})\neq1$ then the estimate for the constant $\gamma$ can be derived from the already proven inequality (\ref{E-formprot}) by applying the latter to the matrix set
\[
\setA'=\rho^{-1}({\setA})\setA=
\{\rho^{-1}({\setA})A_{1},\ldots,\rho^{-1}({\setA})A_{r}\},
\]
for which $\rho({\setA'})=1$. In this case the constant
$\gamma$ in (\ref{E-wirth-irr}) can be estimated from below as follows:
\[
\gamma\ge\frac{\rho^{-(d-1)}({\setA})p_{1}({\setA})\cdots p_{d-1}({\setA})}%
{\left(1+\rho^{-1}({\setA})\|\setA\|\right)^{d-1}}=
\frac{p_{1}({\setA})\cdots p_{d-1}({\setA})}%
{\left(\rho({\setA})+\|\setA\|\right)^{d-1}},
\]
and, if to take into account that by (\ref{Eq-sprad})
$\rho({\setA})\le\|\setA\|$, then
\[
\gamma\ge\frac{p_{1}({\setA})\cdots p_{d-1}({\setA})}%
    {\left(2\|\setA\|\right)^{d-1}}.
\]

\section{Main Theorem}\label{S-irred}

Our aim is to obtain one more explicit a priori estimate of the joint spectral radius with the help of the generalized Gelfand formula. Such an estimate might have many implications in the area of the joint spectral radius. In particular, knowledge of the constant $\gamma$  might be useful in evaluating the local Lipschitz constant of the joint spectral radius at an irreducible inclusion \cite{Wirth:LAA02}.

Denote by ${\setA}_{n}$ for $n\ge 1$ the collection of all finite products of matrices from $\setA\bigcup \{I\}$,
consisting of no more than $k$ factors, that is ${\setA}_{n}=\cup_{k=0}^{n}\setA^{k}$. Then ${\setA}_{n}(x)$
denotes the set of all the vectors $Ax$ where $A\in{\setA}_{n}$. Let $\|\cdot\|$ be a norm in $\mathbb{R}^{d}$; the ball of radius $t$ in this norm is denoted by $\mathbf{S}(t)$.

Let us call the \emph{$p$-measure of irreducibility} of the matrix set $\setA$ (with respect to the norm $\|\cdot\|$) the quantity $\chi_{p}({\setA})$ determined as
\[
\chi_{p}({\setA}) = \inf_{\substack{x\in\mathbb{R}^{d}\\ \|x\|=1}}
\sup \{ t:\mathbf{S}(t)\subseteq{\conv} \{ {\setA}_{p}(x)\cup {\setA}_{p}(-x)\} \}.
\]

Under the name `the measure of quasi-controllability' the measure of irreducibility $\chi_{p}({\setA})$  was introduced and investigated in  \cite{KozPok:DAN92:e,KozPok:TRANS97,KKP:93,KozPok:CADSEM96-005,KKP:CESA98} where the overshooting effects for the transient regimes of  linear remote control systems were studied. The reason why the quantity $\chi_{p}({\setA})$ got the name `the measure of irreducibility' is in the following lemma.

\begin{lemma}\label{L-qcontr} Let $p\ge d-1$. The matrix set $\setA$ is irreducible iff $\chi_{p}({\setA})> 0$.
\end{lemma}

The proof of Lemma~\ref{L-qcontr} can be found in \cite[Thm. 2.4]{KozPok:CADSEM96-005}, \cite[Thm 1.4]{KozPok:TRANS97}, \cite[Thm. 2]{KKP:CESA98}.

\begin{theorem}\label{T-main}
For any $n\ge1$, $p\ge d-1$ the following inequalities are valid
\begin{equation}\label{E-mainineq}
\rho({\setA})\le\|\setA^{n}\|^{1/n}\le \left(\eta_{p}({\setA})\right)^{1/n}\rho({\setA}),
\end{equation}
where
\[
\eta_{p}({\setA})=\frac{\max\{1,\left(\rho({\setA})\right)^{p}\}}
{\chi_{p}({\setA})},
\]
and therefore
\begin{equation}\label{E-compineq}
\left(\nu_{p}({\setA})\right)^{-1/n}
\|\setA^{n}\|^{1/n}\le \rho({\setA})\le \|\setA^{n}\|^{1/n},
\end{equation}
where
\[
\nu_{p}({\setA})=
\frac{\max\{1,\|\setA\|^{p}\}}%
{\chi_{p}({\setA})}.
\]
\end{theorem}

Remark that while computation of $\eta_{p}({\setA})$ requires the knowledge of the value of $\rho({\setA})$, which is a priory unknown, the quantity $\nu_{p}({\setA})$ \textsl{can be evaluated in a finite number of algebraic operations involving only information about the entries of the matrices from $\setA$}.
\bigskip

\section{Proof of Theorem~\protect\ref{T-main}}\label{PRT11}

The proof of Theorem~\ref{T-main} follows the idea of the proof of Theorem 2.3 from \cite{KozPok:TRANS97}, see also \cite[Thm. 2.3]{KozPok:CADSEM96-005}, \cite[Thm. 4]{KKP:CESA98}.

As was noted above the estimate $\rho({\setA})\le\|\setA^{n}\|^{1/n}$ follows from (\ref{E-JSR}). Therefore only the estimate $\|\setA^{n}\|^{1/n}\le \left(\eta_{p}({\setA})\right)^{1/n}\rho({\setA})$ needs to be proved. Suppose the contrary. In this case, for some $n\ge1$,
\[
\|\setA^{n}\|^{1/n}> \left(\eta_{p}({\setA})\right)^{1/n}\rho({\setA}),
\]
or, what is the same,
\[
\|\setA^{n}\|> \eta_{p}({\setA})\rho^{n}({\setA}).
\]
Then, by definition of $\|\setA^{n}\|$, there are matrices $A_{i_{1}},A_{i_{2}},\ldots,A_{i_{n}}\in\setA$ for which
\[
\|A_{i_{n}}\cdots A_{i_{2}}A_{i_{1}}\|=\|\setA^{n}\| > \eta_{p}({\setA})\rho^{n}({\setA}).
\]
Hence, there is a non-zero vector $x_{*}\in\mathbb{R}^{d}$ such that
\[
\|A_{i_{n}}\cdots A_{i_{2}}A_{i_{1}}x_{*}\|> \eta_{p}({\setA})\rho^{n}({\setA})\|x_{*}\|.
\]
It is convenient to rewrite this last inequality in the following form
\begin{equation}\label{E-mu}
\|A_{i_{n}}\cdots A_{i_{2}}A_{i_{1}}x_{*}\|\ge \mu\eta_{p}({\setA})\rho^{n}({\setA})\|x_{*}\|,
\end{equation}
where $\mu$ is some number strictly greater than unity: $\mu>1$.

To complete the proof of Theorem~\ref{T-main} we will need two auxiliary statements. Denote by $\setAinf$ the collection of all finite products of matrices from $\setA$, that is $\setAinf=\cup_{k\ge1}\setA^{k}$.
For each matrix $A\in{\setAinf}$ the amount of the factors $A_{1},A_{2},\ldots,A_{q}\in{\setA}$ in the representation $A = A_{1}A_{2}\cdots A_{q}$ is called the \emph{length of $A$} and is denoted by $\len(A)$.

\begin{lemma}\label{increase}
Given $p\ge d-1$ and the irreducible matrix set $\setA$. Let for some $x_{*}\neq0\in\mathbb{R}^{d}$,
$F\in{\setAinf}$ and $\nu>0$ the inequality
\begin{equation}\label{EL1}
\|Fx_{*}\| \ge \nu\|x_{*}\|,
\end{equation}
hold. Then for any $x\neq0\in\mathbb{R}^{d}$ there is a matrix
$F_x=FL_{x}\in{\setAinf}$ such that $\len(F)\le\len(F_x)\le\len(F)+p$ and
\[
\|F_{x}x\|\ge
\nu\chi_{p}({\setA})\|x\|.
\]
\end{lemma}

\proof Choose a fixed but otherwise arbitrary vector $x\in\mathbb{R}^{d}$. By definition of the measure of irreducibility the vector
$\chi_{p}({\setA})\|x_{*}\|^{-1}x_{*}$ belongs to the centered convex hull of the vectors  ${\setA}_{p}\left(\|x\|^{-1}x\right)$. Then there exist numbers $\theta_{1}$, $\theta_{2}$, \ldots , $\theta_{s}$,
\begin{equation}\label{E127}
\sum^{s}_{i=1}|\theta_{i}|\le 1,
\end{equation}
and matrices
$G_{1},G_{2},\ldots ,G_{s}\in{\setA}_{p}$ such that
\[
\sum^{s}_{i=1}\theta_{i}\|x\|^{-1}G_{i}x
=\chi_{p}({\setA})\|x_{*}\|^{-1}x_{*}.
\]
Therefore
\[
\sum^{s}_{i=1}\theta_{i}\|x\|^{-1}FG_{i}x
=\chi_{p}({\setA})\|x_{*}\|^{-1}Fx_{*},
\]
from which by (\ref{EL1})
\[
\sum^{s}_{i=1}\|\theta_{i}FG_{i}x\|\ge\nu\chi_{p}({\setA})
\|x\|.
\]
Then by (\ref{E127}) there is an index $i$, $1\le i\le s$, for which the matrix $F_{x}=FG_{i}$ satisfies the inequality
\[
\|F_{x}x\|\ge\nu\chi_{p}({\setA})\|x\|.
\]

It remains to note that due to $G_{i}\in{\setA}_{p}$ independently of $x$ the length of the matrix $F_{x}=FG_{i}$ does not exceed the length of the matrix $F$ augmented by
$p$, that is $\len(F)\le\len(F_x)\le\len(F)+p$.  The lemma is proved. \qed

\begin{lemma}\label{increasec}
Let for some set of matrices $A_{i_{1}},A_{i_{2}},\ldots,A_{i_{n}}\in\setA$ the inequality (\ref{E-mu}) hold. Then for any  $x\neq0\in\mathbb{R}^{d}$ there is a matrix $F_x\in{\setAinf}$ such that $n\le\len(F_x)\le n+p$ and
\begin{equation}\label{E-lenestc}
\|F_{x}x\|\ge\left(\eta\rho({\setA})\right)^{\len(F_{x})}\|x\|,
\end{equation}
where $\eta=\mu^{1/(n+p)}>1$.
\end{lemma}

\proof Set $\nu=\mu\eta_{p}({\setA})\rho^{n}({\setA})$. Then by Lemma~\ref{increase} for any  $x\neq0\in\mathbb{R}^{d}$ there ia a matrix
$F_x\in{\setAinf}$ such that $n\le\len(F_x)\le n+p$ and
\[
\|F_{x}x\|\ge\mu\eta_{p}({\setA})\chi_{p}({\setA})\rho^{n}({\setA})\|x\|,
\]
which is equivalent, by definition of $\eta_{p}({\setA})$, to the following chain of relations
\begin{multline}\label{E-fin}
\|F_{x}x\|\ge\mu\frac{\max\{1,\rho^{p}({\setA})\}}
{\chi_{p}({\setA})}\chi_{p}({\setA})\rho^{n}({\setA})\|x\|=\\
\mu\max\{1,\rho^{p}({\setA})\}\rho^{n}({\setA})\|x\|=
\frac{\mu}{\eta^{\len(F_{x})}}\cdot \frac{\max\{1,\rho^{p}({\setA})\}}{\rho^{\len(F_{x})-n}({\setA})}\cdot \left(\eta\rho({\setA})\right)^{\len(F_{x})}=\\
\mu^{1-\frac{\len(F_{x})}{n+p}}\cdot \frac{\max\{1,\rho^{p}({\setA})\}}{\rho^{\len(F_{x})-n}({\setA})}\cdot \left(\eta\rho({\setA})\right)^{\len(F_{x})}.
\end{multline}

As is easy to see, by the inequalities $n\le\len(F_x)\le n+p$ as the first as the second factors in the right-hand part of (\ref{E-fin}) are greater than or equal to $1$, from which the required inequality (\ref{E-lenestc}) follows. \qed

\begin{lemma}\label{expo} Let for some set of matrices $A_{i_{1}},A_{i_{2}},\ldots,A_{i_{n}}\in\setA$ the inequality (\ref{E-mu}) hold. Then there is a sequence of matrices  $H_{k}\in\setA$, $k=0,1,\ldots,$ for which  $\len(H_{k})\to\infty$ as $k\to\infty$ and
\begin{equation}\label{E-Hnorms}
\|H_{k}\|\ge \left(\eta\rho({\setA})\right)^{\len(H_{k})},\qquad k=0,1,\ldots~.
\end{equation}
\end{lemma}

\proof Let us choose a non-zero but otherwise arbitrary vector $x\in\mathbb{R}^{d}$, and construct an auxiliary sequence of the vectors $\{x(k)\}$, $k=0,1,\ldots $, by setting $x(0)=x$,
\[
x(k+1)=F_{x(k)}x(k),\qquad k=0,1,\ldots ,
\]
where $F_{x}$ is the matrix defined in Lemma~\ref{increasec}. Then by Lemma~\ref{increasec} the following relations hold:
\[
\|x(k+1)\|=\|F_{x(k)}x(k)\| \ge\left(\eta\rho({\setA})\right)^{\len(F_{x(k)})}\|x(k)\|,\qquad k=0,1,\ldots~,
\]
from which
\begin{multline}\label{E-finest}
\|F_{x(k)}\cdots F_{x(1)}F_{x(0)}x(0)\| \ge\\%
\left(\eta\rho({\setA})\right)^{\len(F_{x(k)})+\dots+\len(F_{x(1)})+\len(F_{x(0)})}\|x(0)\|,\qquad k=0,1,\ldots~,
\end{multline}
By setting now
\[
H_{k}=F_{x(k)}\cdots F_{x(1)}F_{x(0)},\qquad k=0,1,\ldots~,
\]
from (\ref{E-finest}) we get
\[
\|H_{k}x(0)\| \ge\left(\eta\rho({\setA})\right)^{\len(H_{k})}\|x(0)\|,\qquad k=0,1,\ldots~,
\]
which imply the required estimates (\ref{E-Hnorms}). The lemma is proved. \qed

We now complete the proof of Theorem~\ref{T-main}. By Lemma~\ref{expo} the inequality (\ref{E-mu}) implies the existence of such a sequence of matrices $H_{k}\in\setA$, $k=0,1,\ldots$, for which the estimates (\ref{E-Hnorms}) hold. Then, by setting $n_{k}=\len(H_{k})$ from the definition of $\|\setA^{n}\|$ we get
\[
\|\setA^{n_{k}}\|^{1/n_{k}}\ge \|H_{k}\|^{1/n_{k}}\ge\eta\rho({\setA}).
\]
Passing here to the limit in the left-hand part as $n_{k}\to\infty$ we obtain:
\[
\rho({\setA})=\limsup_{n\to\infty}\|\setA^{n}\|^{1/n} \ge \limsup_{n_{k}\to\infty}\|\setA^{n_{k}}\|^{1/n_{k}}\ge \eta\rho({\setA}),
\]
which is impossible since $\eta>1$ and the value of $\rho({\setA})$ is strictly positive for any irreducible matrix set (the latter follows, for example, from the inequality (\ref{E-extnorm}) valid for some extremal norm). The obtained contradiction completes the proof of the estimates (\ref{E-mainineq}).

It remains only to note that the estimates (\ref{E-compineq}) directly follow from (\ref{E-mainineq}) and the inequality $\rho({\setA})\le\|\setA\|$ valid by (\ref{Eq-sprad}).

The proof of Theorem~\ref{T-main} is completed.

\section{Comments and Examples}\label{S-ex}

Remark that with one trivial exception Theorem~\ref{T-main} is inapplicable for matrix sets consisting of a single matrix. This is because a singleton matrix set $\setA=\{A\}$ may be irreducible only in the case when $A$ is a matrix of dimension $2\times 2$, and also it has no real eigenvalues.

In a sense, irreducibility of matrix sets naturally arising in control theory is similar to the Kalman `complete controllability plus complete observability' property.

Consider, in particular, a real $d\times d$ matrix $A=(a_{ij})$ and column-vectors $b,c\in\mathbb{R}^{d}$. Form with their help the matrix set
${\setA}={\setA}(A,b,c)=\{A,Q\}$ where
$Q=(q_{ij})=bc^{T}$ is the matrix with the entries  $q_{ij}=b_{i}c_{j}$,
$i,j=1,\ldots,d$.
Then the matrix set ${\setA}(A,b,c)$ is irreducible if and only if the pair $(A,b)$ is completely controllable and the pair $(A,c)$
is completely observable by Kalman \cite{KozPok:TRANS97}.

Reproduce a pair of examples from \cite{KozPok:TRANS97} which demonstrate how to verify the irreducibility of a matrix set and how to evaluate its measure of irreducibility. Proofs of the corresponding statements can be found in \cite{KozPok:TRANS97}. The next example came from the theory of asynchronous systems, see, e.g., \cite{AKKK:92:e,KKKK:MCS84}. Given a real $d\times d$ matrix
$A= (a_{ij})$, form the matrix set
${\setP}(A) = \{A_{1}, A_{2}, \ldots , A_{d} \}$ where
\[
A_{i} =\left(\begin{array}{cccccc}
1 & 0 &\dots & 0 &\dots & 0\\ 0 &
$1$ &\dots & 0 &\dots & 0\\
\vdots &\vdots &\ddots &\vdots &\ddots &\vdots\\
a_{i1}& a_{i2}&\dots & a_{ii}&\dots & a_{id}\\
\vdots &\vdots &\ddots &\vdots &\ddots &\vdots\\
0 & 0 &\dots & 0 &\dots & 1\end{array}\right).
\]

Recall that the matrix $A$ is said to be \emph{indecomposable} if it cannot be represented in a block-triangle form
\[
A=\left(\begin{array}{cc} B&C\\0&D\end{array}\right)
\]
by any transposition of its rows and corresponding columns.
The matrix $A$ is indecomposable if and only if it has no proper non-zero invariant set spanned over some number of the basis vectors
\[
e_{i} = \{0,0, \ldots,1, \ldots , 0 \},\qquad i = 1,2,\ldots ,d.
\]

Let $\|\cdot\|$ be the norm in $\mathbb{R}^{d}$ defined by the equality $\|x\|
=|x_{1}| + |x_{2}| +\ldots + |x_{d}|$. Set
\[
\alpha ={\frac{1}{2d}}\min \{\|(A-I)x\|:\|x\|=1 \},\qquad
\beta ={\frac{1}{2}}\min \{|a_{ij}| : i\neq j, a_{ij}\neq 0 \}.
\]
\begin{example}\label{Ex11} The matrix set ${\setP}(A)$ is irreducible iff the matrix $A$ indecomposable and
$1$ is not its eigenvalue. In this case
${\chi}_{d}[{\setP}(A)]\ge {\alpha}{\beta}^{d-1}$.
\end{example}

The next example was communicated to the author by E. Kaszkurewicz. Let ${\setV}(A)$ be the matrix set composed of the matrices $A$,
$V_{1}A$, $V_{2}A$,\ldots , $V_{d}A$. Here $A$ is a real
$d\times d$ matrix with the entries $a_{ij}$, and $V_{i}$, $i=1,2,\ldots ,d$, are the diagonal matrices of the form
\[
V_{i}= \diag \{v_{1i},v_{2i},\ldots ,v_{ii},\ldots ,v_{di} \},
\]
where $v_{ij}=1$ for $i\neq j$, $v_{ij}=-1$ for $i=j$.

Let again $\|\cdot\|$ be the norm in $\mathbb{R}^{d}$ defined by the equality $\|x\|
=|x_{1}| + |x_{2}| +\ldots + |x_{d}|$. Set
\[
\tilde\alpha ={\frac{1}{d}}\min \{\|Ax\|:\ \|x\|=1 \},\qquad
\tilde\beta =\min \{|a_{ij}| :\  i\neq j, a_{ij}\neq 0 \}.
\]

\begin{example}\label{EX211} The matrix set ${\setV}(A)$
is irreducible iff the matrix $A$ is indecomposable and non-degenerate. In this case ${\chi}_{d}[{\setV}(A)]\ge{\tilde\alpha}{\tilde\beta}^{d-1}$.
\end{example}

\section{Acknowledgments}

The author is grateful to the reviewer for a number of valuable remarks.

\end{document}